\documentclass[11pt,a4paper,reqno]{amsart}
\usepackage[T1]{fontenc}
\usepackage{lmodern}
\usepackage{amsmath,amssymb,mathtools}
\usepackage[a4paper,margin=29mm]{geometry}
\usepackage{microtype}
\usepackage{xcolor}
\definecolor{linkblue}{RGB}{31,63,92}
\usepackage[colorlinks=true,allcolors=linkblue]{hyperref}
\hypersetup{pdftitle={C*(F2 x F2) fails the local lifting property},pdfsubject={Local lifting and uniform intertwining obstructions}}
\numberwithin{equation}{section}
\newtheorem{theorem}{Theorem}[section]
\newtheorem{maintheorem}{Theorem}
\newtheorem{maincorollary}[maintheorem]{Corollary}

\newtheorem{proposition}[theorem]{Proposition}
\newtheorem{lemma}[theorem]{Lemma}

\theoremstyle{remark}
\newtheorem*{remark}{Remark}
\newcommand{\C}{\mathbb C}
\newcommand{\F}{\mathbb F}
\newcommand{\Z}{\mathbb Z}
\newcommand{\B}{\mathbb B}
\newcommand{\A}{\mathcal A}
\newcommand{\id}{\mathrm{id}}
\newcommand{\tmax}{\otimes_{\max}}
\DeclareMathOperator{\tr}{tr}
\DeclareMathOperator{\Tr}{Tr}
\DeclareMathOperator{\Ad}{Ad}
\DeclareMathOperator{\Ree}{Re}
\DeclareMathOperator{\SL}{SL}
\title[\(C^*(\F_2\times\F_2)\) fails the local lifting property]{\texorpdfstring{\(C^*(\F_2\times\F_2)\)}{C*(F2 x F2)} fails the local lifting property}
\date{November 17, 2026}
\begin{document}
\author{Jean ROYDOR}
\address{Institut de Math\'ematiques de Jussieu, Sorbonne Universit\'e}
\email{jroydor@imj-prg.fr}
\begin{abstract}
We prove that the full group \(C^*\)-algebra \(C^*(\F_k\times\F_m)\)
fails the local lifting property for all integers \(k,m\ge2\), and
therefore disprove Ozawa's Strong Kirchberg Conjecture. Our proof is
inspired by a construction of Ioana concerning Hilbert--Schmidt stability.
We then strengthen this result by showing that \(C^*(\F_k\times\F_m)\)
does not locally embed into \(C^*(\F_\infty)\). Consequently, the lifting
property, the local lifting property, and local embeddability into
\(C^*(\F_\infty)\) are not preserved under maximal tensor products.
\end{abstract}

\maketitle

\section{Introduction}
Kirchberg introduced the local lifting property in \cite{Kirchberg}.
A unital \(C^*\)-algebra \(A\) has the \emph{local lifting property}
(LLP) if every unital completely positive (u.c.p.) map
\(\alpha:A\to C/J\), where \(J\) is a closed ideal in a unital
\(C^*\)-algebra \(C\), admits a u.c.p.\ lift on each finite-dimensional
operator system \(E\subset A\). The \emph{lifting property} (LP)
requires a u.c.p.\ lift on all of \(A\).

The duality between the LLP and the weak expectation property played a
central role in Connes' embedding problem; see \cite{Pisier}.
Junge and Pisier proved that \(\B(\ell_2)\) fails the LLP \cite{JP}.
Amenable and free groups have full group \(C^*\)-algebras with the LLP.
Thom gave the first explicit group whose full \(C^*\)-algebra fails
the LLP; his group has property \((T)\) \cite{Thom}. Ioana, Spaas and Wiersma combined relative
property \((T)\) with cohomological obstructions to prove that further
full group \(C^*\)-algebras, including \(C^*(\SL_3(\Z))\), fail the LLP
\cite[Corollary~B]{ISW}.

In 2004, Ozawa asked whether maximal tensor products preserve the LLP
and observed that it suffices to consider \(C^*(\F_\infty\times\F_\infty)\)
\cite[end of Section~3]{Ozawa04}. He later reduced the question to
\(C^*(\F_2\times\F_2)\) and related it to his Strong Kirchberg
Conjecture \cite[discussion after Theorem~4]{Ozawa13}. Enders and Shulman
proved that this algebra has the LLP if and only if
\(\operatorname{Ext}(C^*(\F_2\times\F_2))\) is a group
\cite[Corollary~7.6]{ES}. The question also appears in
\cite[Question~1.6]{FFW} and among the open problems in
Pisier's book \cite[Chapter~25]{Pisier}.

\begin{maintheorem}\label{thm:main}
For integers \(k,m\ge2\), the full group \(C^*\)-algebra
\(C^*(\F_k\times\F_m)\) fails the LLP.
\end{maintheorem}

Since \(\F_2\times\F_2\) has the Haagerup property, it has no infinite
subgroup \(H\) for which \((\F_2\times\F_2,H)\) has relative property
\((T)\). To our knowledge, this is the first example of a group without
such a subgroup whose full group \(C^*\)-algebra fails the LLP.
Theorem~\ref{thm:main} also answers \cite[Question~1.8]{FFW}.

For Hilbert spaces \(H\subset\widetilde H\), let
\(\Phi_H:\B(\widetilde H)\to\B(H)\) be the compression
\(\Phi_H(X)=P_HX|_H\), where \(P_H\) is the orthogonal projection onto
\(H\). Ozawa's reformulation of the LLP \cite[Lemma~5]{Ozawa13}
gives the following consequence.

\begin{maincorollary}\label{cor:ozawa-matricial}
There exists \(\kappa>0\) such that, for every integer \(m\ge2\) and
\(\varepsilon>0\), there are an integer \(N\ge1\) and unitaries
\(U_1,\ldots,U_m,V_1,\ldots,V_m\in\mathcal U(\ell_2^N)\) with
\[
 \|[U_i,V_j]\|\le\varepsilon\qquad(1\le i,j\le m),
\]
such that, for any Hilbert space \(\widetilde H\supset\ell_2^N\) and unitaries
\(\widetilde U_1,\ldots,\widetilde U_m,\allowbreak
  \widetilde V_1,\ldots,\widetilde V_m\in\mathcal U(\widetilde H)\)
satisfying
\[
 [\widetilde U_i,\widetilde V_j]=0\qquad(1\le i,j\le m),
\]
one has
\[
 \max_{1\le i\le m}
 \bigl\{\|\Phi_{\ell_2^N}(\widetilde U_i)-U_i\|,
         \|\Phi_{\ell_2^N}(\widetilde V_i)-V_i\|\bigr\}>\kappa.
\]
\end{maincorollary}

A \(C^*\)-algebra \(A\) \emph{locally embeds into \(C^*(\F_\infty)\)}
if, for every finite-dimensional subspace \(E\subset A\) and every
\(\varepsilon>0\), there is an injective linear map
\(j:E\to C^*(\F_\infty)\) with
\[
 \|j\|_{\mathrm{cb}}\,\|j^{-1}:j(E)\to E\|_{\mathrm{cb}}<1+\varepsilon.
\]
Here \(\|\cdot\|_{\mathrm{cb}}\) is the completely bounded norm.
The LLP implies local embeddability into \(C^*(\F_\infty)\).
As explained in \cite{Roydor}, local embeddability is an analogue of
QWEP on the LLP side of the LLP/WEP duality. We obtain the following strengthening.
\enlargethispage{2\baselineskip}

\begin{maintheorem}\label{thm:local-embedding}
For integers \(k,m\ge2\), the full group \(C^*\)-algebra
\(C^*(\F_k\times\F_m)\) does not locally embed into \(C^*(\F_\infty)\).
\end{maintheorem}

Each \(C^*(\F_k)\) has the LP and locally embeds into \(C^*(\F_\infty)\).
Since \(C^*(\F_k\times\F_m)=C^*(\F_k)\tmax C^*(\F_m)\),
Theorems~\ref{thm:main} and~\ref{thm:local-embedding} show that none of
these three properties is preserved by maximal tensor products.

Our construction is inspired by Ioana's work on Hilbert--Schmidt stability.
He constructed a \(*\)-homomorphism from \(C^*(\F_2\times\F_2)\) into a
tracial ultraproduct of matrix algebras with no \(*\)-homomorphic lift,
and asked whether it admits a u.c.p.\ lift \cite[p.~4031]{Ioana}.
We use three generators in each factor: two provide spectral gap estimates,
while the third encodes the matrix families associated with a deformation
of Walsh bases. The index-two inclusion
\(\F_3\hookrightarrow\F_2\) then gives the rank-two case.

Let \(u_1,u_2,u_3\) be the free generators of \(\F_3\), and let \(e\)
be its identity.

\begin{maintheorem}\label{thm:gap}
There are \(c>0\) and finite-dimensional Hilbert spaces
\((H_n)_{n\in\mathbb N}\), each carrying two triples
\(U_{i,n},V_{i,n}\in\mathcal U(H_n)\), \(i=1,2,3\), such that
\[
 \max_{1\le i,j\le3}\|[U_{i,n},V_{j,n}]\|\longrightarrow0,
\]
and, for every \(n\), every unitary representation
\(\pi:\F_3\times\F_3\to\mathcal U(K)\) on a Hilbert space \(K\),
and every \(T\in\B(H_n,K)\),
\begin{equation}\label{eq:gap}
 \sum_{i=1}^3\bigl(
 \|\pi(u_i,e)T-TU_{i,n}\|_{\mathrm{HS}}^2
 +\|\pi(e,u_i)T-TV_{i,n}\|_{\mathrm{HS}}^2\bigr)
 \ge c\|T\|_{\mathrm{HS}}^2.
\end{equation}
\end{maintheorem}
Here \(\|T\|_{\mathrm{HS}}^2=\Tr(T^*T)\), with \(\Tr\) unnormalized.

Let
\[
 q:\ell_\infty(\B(H_n))\longrightarrow
     \ell_\infty(\B(H_n))/c_0(\B(H_n))
\]
be the quotient map. The asymptotic commutation relations give a unitary
representation \(\rho\) of \(\F_3\times\F_3\) in this quotient, defined by
\[
 \rho(u_i,e)=q((U_{i,n})_n),\qquad
 \rho(e,u_i)=q((V_{i,n})_n)\qquad(i=1,2,3).
\]

%We identify group elements with their canonical unitaries in full group \(C^*\)-algebras.

The two triples are adapted from Ioana's block construction
\cite[proof of Proposition~1.5, Section~5.1]{Ioana}, whose commutators
already tend to zero in operator norm
\cite[Lemma~6.3 and Remark~1.9]{Ioana}. We use pairs with spectral gap
for conjugation on both outer tensor factors \cite[Lemma~4.3]{Ioana}.

Like Ioana's deformation, our construction attenuates tensor components
according to their degree
\cite[Construction~1.7, Lemma~6.2 and Eq.~(6.1)]{Ioana}.
We realize this mechanism with two Walsh bases of maximal abelian
subalgebras in one matrix algebra, obtained by rotating each tensor
factor; Ioana works in a doubled tensor product. He uses perturbation
results and incompatible polynomial and exponential dimension bounds
\cite[Corollary~6.4 and Lemmas~6.5--6.6]{Ioana}.
We replace this dimension argument by an estimate on a single vector
associated with an almost intertwiner through standard form
(Lemma~\ref{lem:standard} and Proposition~\ref{prop:quant-reduction}).
The deformation controls its Walsh degree
(Proposition~\ref{prop:walsh-degree-control}); the matrix spectral gap
and an exact exchange of complementary subspaces yield uniform
coercivity (Lemma~\ref{lem:quant-complementary} and
Proposition~\ref{prop:quant-walsh-control}), giving
\eqref{eq:quant-concrete-coercivity} on the associated Walsh space.
This yields the intertwining obstruction \eqref{eq:gap} of
Theorem~\ref{thm:gap} for representations on arbitrary Hilbert spaces.

\section{The construction of the six unitaries}\label{sec:matrices}

\subsection{Two Walsh bases on one cube}
We write \(\tr_k=k^{-1}\Tr\) for the normalized trace on
\(\mathbb M_k\). For a unital sub-\(C^*\)-algebra
\(D\subset\mathbb M_k\), let \(L_2(D)\) denote \(D\) with inner
product \(\langle b,c\rangle=\tr_k(c^*b)\) and norm \(\|b\|_2\).
For \(n\ge2\), put \([n]=\{1,\ldots,n\}\) and
\(H_n^0=(\ell_2^2)^{\otimes n}\), so that
\(\B(H_n^0)=\mathbb M_2^{\otimes n}\cong\mathbb M_{2^n}\), and set
\[
 \sigma=\begin{pmatrix}1&0\\0&-1\end{pmatrix},\qquad
 \tau=\begin{pmatrix}0&1\\1&0\end{pmatrix},\qquad
 \sigma_n=(1-1/n)^{1/2}\sigma+n^{-1/2}\tau\in\mathbb M_2.
\]
Since \(\sigma\tau=-\tau\sigma\), \(\sigma_n\) is a traceless
self-adjoint unitary. For \(k\in[n]\), define
\[
\begin{aligned}
 x_{k,n}&=1_2^{\otimes(k-1)}\otimes\sigma\otimes1_2^{\otimes(n-k)}
          \in\mathbb M_{2^n},\\
 y_{k,n}&=1_2^{\otimes(k-1)}\otimes\sigma_n\otimes1_2^{\otimes(n-k)}
          \in\mathbb M_{2^n},
\end{aligned}
\]
where \(1_2\) is the identity matrix and factors with exponent zero
are omitted. For \(S\subseteq[n]\), put
\[
 \chi_{S,n}=\prod_{k\in S}x_{k,n},\qquad
 \psi_{S,n}=\prod_{k\in S}y_{k,n}\quad\text{in }\mathbb M_{2^n},
 \qquad \chi_{\varnothing,n}=\psi_{\varnothing,n}=1_{H_n^0}.
\]
The groups
\[
 G_n=\{\chi_{S,n}:S\subseteq[n]\},\qquad
 F_n=\{\psi_{S,n}:S\subseteq[n]\}\subset\mathcal U(H_n^0)
\]
consist of real matrices. Their spans
\(\A_n=\operatorname{span}G_n\) and
\(\mathcal B_n=\operatorname{span}F_n\) are unital abelian
sub-\(C^*\)-algebras of \(\mathbb M_{2^n}\).
The groups \(G_n\) and \(F_n\) form orthonormal bases of
\(L_2(\A_n)\) and \(L_2(\mathcal B_n)\), respectively, called their
\emph{Walsh bases}. We assign degree \(|S|\) to \(\chi_{S,n}\) and
\(\psi_{S,n}\). Their commutants satisfy
\begin{gather*}
 G_n'\cap\mathbb M_{2^n}=\A_n,\\
 F_n'\cap\mathbb M_{2^n}=\mathcal B_n.
\end{gather*}
Writing \(\triangle\) for symmetric difference, we have
\begin{equation}\label{eq:toggle}
 x_{k,n}\chi_{S,n}=\chi_{S\triangle\{k\},n}
 \qquad\text{for every }k\in[n]\text{ and }S\subseteq[n].
\end{equation}

\begin{lemma}\label{lem:cube}
For \(k\in[n]\) and \(z\in G_n\cup F_n\),
\(\|[x_{k,n},z]\|\le2n^{-1/2}\).
The trace-preserving conditional expectation
\(E_{\mathcal B_n}:\mathbb M_{2^n}\to\mathcal B_n\) satisfies
\begin{equation}\label{eq:projection}
 E_{\mathcal B_n}(\chi_{S,n})=(1-1/n)^{|S|/2}\psi_{S,n}
 \qquad\text{for every }S\subseteq[n].
\end{equation}
\end{lemma}
\begin{proof}
The only nonzero commutators occur for \(z=\psi_{T,n}\) with \(k\in T\),
and have norm \(\|[\sigma,\sigma_n]\|=2n^{-1/2}\).
The tensor product trace gives
\[
 \langle\chi_{S,n},\psi_{T,n}\rangle
 =\begin{cases}
 (1-1/n)^{|S|/2},&S=T,\\
 0,&S\ne T.
 \end{cases}
\]
Since \(E_{\mathcal B_n}\) is the orthogonal projection onto
\(L_2(\mathcal B_n)\), this proves \eqref{eq:projection}.
\end{proof}

For a vector in \(L_2(\mathcal A_n)\), proximity to \(L_2(\mathcal B_n)\) controls its Walsh degree
(including for Hilbert space coefficients), this will help us construct an almost invariant vector in the Walsh space with small Walsh number.

\begin{proposition}\label{prop:walsh-degree-control}
Let \(H\) be a Hilbert space, let \(n\ge2\), and write
\[
 x=\sum_{S\subseteq[n]}c_S\otimes\chi_{S,n}
   \in H\otimes L_2(\A_n),\qquad c_S\in H.
\]
For every \(\delta>0\),
\begin{equation}\label{eq:walsh-tail-control}
 \sum_{\substack{S\subseteq[n]\\|S|\ge\delta n}}\|c_S\|^2
 \le
 (1-e^{-\delta})^{-1}
 \operatorname{dist}\bigl(x,H\otimes L_2(\mathcal B_n)\bigr)^2.
\end{equation}
Moreover,
\begin{equation}\label{eq:walsh-degree-control}
 \sum_{S\subseteq[n]}\frac{|S|}{n}\|c_S\|^2
 \le
 (1-e^{-1})^{-1}
 \operatorname{dist}\bigl(x,H\otimes L_2(\mathcal B_n)\bigr)^2.
\end{equation}
The distances are taken in \(H\otimes L_2(\mathbb M_{2^n})\).
\end{proposition}
\begin{proof}
By \eqref{eq:projection} and orthogonality,
\begin{equation}\label{eq:walsh-distance}
 \begin{aligned}
 \operatorname{dist}\bigl(x,H\otimes L_2(\mathcal B_n)\bigr)^2
 &=\|x\|^2-\|(1_H\otimes E_{\mathcal B_n})x\|^2\\
 &=\sum_{S\subseteq[n]}
      \bigl[1-(1-1/n)^{|S|}\bigr]\|c_S\|^2.
 \end{aligned}
\end{equation}
For \(0\le p\le n\),
\[
 1-(1-1/n)^p\ge1-e^{-p/n}\ge(1-e^{-1})p/n,
\]
by concavity of \(r\mapsto1-e^{-r}\) on \([0,1]\).
For \(p\ge\delta n\), the same expression is at least
\(1-e^{-\delta}\). This proves both inequalities.
\end{proof}

\subsection{A uniform spectral gap estimate}
By \cite[Lemma~4.3]{Ioana}, there are a constant \(C>0\) and
unitaries
\[
 R_{1,n},R_{2,n}\in\mathcal U(\ell_2^n)\qquad(n\ge2)
\]
such that
\begin{equation}\label{eq:quantumgap}
 \|x-\tr_n(x)1_n\|_2^2
 \le C\sum_{i=1}^2\|R_{i,n}xR_{i,n}^*-x\|_2^2
 \qquad\text{for every }x\in\mathbb M_n.
\end{equation}
Here \(1_n\) is the identity matrix. The estimate also
holds with the unnormalized Hilbert--Schmidt norm and with the unitaries
replaced by their entrywise conjugates. A bar will denote entrywise
complex conjugation in the standard basis.
The map \(x\mapsto\tr_n(x)1_n\) is the orthogonal projection onto
\(\C1_n\). Tensoring the operators in this inequality with the
identity of any Hilbert space preserves the inequality, as follows by
expanding in an orthonormal basis of that space.

\subsection{The two triples}
For \(n\ge2\), put \(m_n=2^{n+1}\) and enumerate first \(G_n\) and
then \(F_n\) as \((z_{l,n})_{l=1}^{m_n}\subset\mathcal U(H_n^0)\),
retaining both occurrences of the identity. Set
\[
 H_n=\ell_2^n\otimes H_n^0\otimes\ell_2^{m_n}.
\]
Using the standard diagonal matrix units
\(e_{kk}\in\mathbb M_n\), \(k\in[n]\), and
\(e_{ll}\in\mathbb M_{m_n}\), \(l\in[m_n]\), define
\(U_{i,n},V_{i,n}\in\mathcal U(H_n)\), \(i=1,2,3\), by
\begin{equation}\label{eq:encoding}
\begin{array}{c@{\qquad}c@{\qquad}c}
 i & U_{i,n} & V_{i,n} \\[2pt]
 \hline
 1 & R_{1,n}\otimes 1\otimes 1 & 1\otimes 1\otimes R_{1,m_n} \\
 2 & R_{2,n}\otimes 1\otimes 1 & 1\otimes 1\otimes R_{2,m_n} \\[2pt]
 3 & \displaystyle\sum_{k=1}^n e_{kk}\otimes x_{k,n}\otimes 1
   & \displaystyle\sum_{l=1}^{m_n} 1\otimes z_{l,n}\otimes e_{ll}.
\end{array}
\end{equation}
Each \(1\) is the identity on the corresponding tensor factor.
All commutators \([U_{i,n},V_{j,n}]\) vanish except possibly
\([U_{3,n},V_{3,n}]\), whose \((k,l)\)-block is
\([x_{k,n},z_{l,n}]\). Lemma~\ref{lem:cube} gives
\begin{equation}\label{eq:smallcomm}
 \max_{1\le i,j\le3}\|[U_{i,n},V_{j,n}]\|
 \le2n^{-1/2}\longrightarrow0.
\end{equation}

The next proposition relates conjugation by \(V_{3,n}\) to the orthogonal
projections associated with \(\mathcal A_n\) and \(\mathcal B_n\), with
Hilbert space coefficients.

\begin{proposition}\label{prop:walsh-averaging}
Let \(H\) be a Hilbert space and \(n\ge2\). Let
\(E_A,E_B\in\B(H\otimes L_2(\B(H_n)))\) be the orthogonal
projections onto
\[
 \begin{aligned}
 &H\otimes L_2(\mathbb M_n)\otimes L_2(\A_n)\otimes\C1_{m_n},\\
 &H\otimes L_2(\mathbb M_n)\otimes L_2(\mathcal B_n)\otimes\C1_{m_n},
 \end{aligned}
\]
respectively. Put \(S_3=1_H\otimes\Ad V_{3,n}
\in\mathcal U(H\otimes L_2(\B(H_n)))\), where
\(\Ad V_{3,n}(x)=V_{3,n}xV_{3,n}^*\) for \(x\in\B(H_n)\).
Then, for every \(h\in H\otimes L_2(\mathbb M_n)
\otimes L_2(\mathbb M_{2^n})\otimes\C1_{m_n}\),
\begin{equation}\label{eq:walsh-averaging}
 \|h-E_Ah\|^2+\|h-E_Bh\|^2=\|S_3h-h\|^2.
\end{equation}
\end{proposition}
\begin{proof}
For \(h\) in the indicated subspace and \((G,E)=(G_n,E_A)\) or
\((F_n,E_B)\), the commutant identities in Section~2.1 give
\[
 Eh=|G|^{-1}\sum_{g\in G}(1_H\otimes1\otimes\Ad g\otimes1)h.
\]
Here \(\Ad g(x)=gxg^*\) for \(x\in\mathbb M_{2^n}\).
Indeed, averaging a unitary representation of a finite group gives
the orthogonal projection onto its fixed vectors. Hence
\[
 \begin{aligned}
 |G|^{-1}\sum_{g\in G}
 \|(1_H\otimes1\otimes\Ad g\otimes1)h-h\|^2
 &=2\|h\|^2-2\Ree\langle Eh,h\rangle\\
 &=2\|h-Eh\|^2.
 \end{aligned}
\]
Since \(|G_n|=|F_n|=m_n/2\), the block formula \eqref{eq:encoding},
with \(\|e_{ll}\|_2^2=m_n^{-1}\) for every \(l\in[m_n]\), gives
\[
 \begin{aligned}
 \|h-E_Ah\|^2+\|h-E_Bh\|^2
 &=m_n^{-1}\sum_{l=1}^{m_n}
    \|(1_H\otimes1\otimes\Ad z_{l,n}\otimes1)h-h\|^2\\
 &=\|S_3h-h\|^2.
 \end{aligned}
\]
\end{proof}

\section{Passing to the standard form}\label{sec:standard}
In the proof of Proposition~\ref{prop:quant-reduction}, we pass to standard
form to eliminate the action of \(\pi\) on the second factor from the
estimates. The following lemma can be viewed as an absorption principle.
For a von Neumann algebra \(M\), let \(M_*^+\) be its cone of
positive normal functionals. Its standard Hilbert space \(L_2(M)\)
carries the left representation \(\lambda_M:M\to\B(L_2(M))\) and
the positive cone \(\mathfrak P_M\subset L_2(M)\).
Each \(\varphi\in M_*^+\) has a unique representative
\(\xi_\varphi\in\mathfrak P_M\) satisfying
\[
 \varphi(x)=\langle\lambda_M(x)\xi_\varphi,\xi_\varphi\rangle
 \qquad\text{for every }x\in M.
\]
Inner products are linear in the first variable.
For a normal unital representation \(\theta:M\to\B(H)\) on a Hilbert
space \(H\) and unit vectors \(v_1,v_2\in H\), let
\(\varphi_j(x)=\langle\theta(x)v_j,v_j\rangle\), \(j=1,2\).
The Araki--Yamagami inequality \cite[Theorem~2]{AY} gives
\begin{equation}\label{eq:AY}
 \|\xi_{\varphi_1}-\xi_{\varphi_2}\|_{L_2(M)}^2
 \le2\|v_1-v_2\|^2.
\end{equation}

The positive cone of \(L_2(\mathbb M_n)\) is \(\mathbb M_n^+\).
For \(z\in\mathbb M_n\) and \(w\in\mathcal U(\mathbb M_n)\), define
\(L_z,\Ad w\in\B(L_2(\mathbb M_n))\) by
\[
 L_zX=zX,\qquad (\Ad w)X=wXw^*\qquad(X\in\mathbb M_n).
\]
The operator \(\Ad w\) is unitary. For
\(N=M\bar\otimes\mathbb M_n\), we use the standard identification
\[
 L_2(N)=L_2(M)\otimes L_2(\mathbb M_n),\qquad
 \lambda_N(x\otimes z)=\lambda_M(x)\otimes L_z.
\]

\begin{samepage}
\begin{lemma}\label{lem:standard}
Let \(M\subset\B(K)\) be a unital von Neumann algebra on a Hilbert
space \(K\), and let \(n\ge1\).
For a unit vector \(v\in K\otimes\ell_2^n\), let
\(\varphi\in N_*^+\) be its vector state on
\(N=M\bar\otimes\mathbb M_n\), and let
\(\xi=\xi_\varphi\in\mathfrak P_N\subset L_2(N)\).
Then, for every \(a\in\mathcal U(M)\), \(b\in\mathcal U(M')\), and
\(z,w\in\mathcal U(\mathbb M_n)\),
\begin{align}
 \|(\lambda_M(a)\otimes L_z)\xi-\xi\|
      &=\|(a\otimes z)v-v\|, \label{eq:leftstandard}\\
 \|(1_{L_2(M)}\otimes\Ad w)\xi-\xi\|^2
      &\le2\|(b\otimes w)v-v\|^2. \label{eq:rightstandard}
\end{align}
\end{lemma}
\end{samepage}
\begin{proof}
Since \(v\) and \(\xi\) represent the same state,
\[
 \|(a\otimes z)v-v\|^2
 =2-2\Ree\,\varphi(a\otimes z)
 =\|(\lambda_M(a)\otimes L_z)\xi-\xi\|^2.
\]
Since \(b\in M'\), the state of \((b\otimes w)v\) on \(N\) is
\(\varphi\circ(\id_M\otimes\Ad(w^*))\). Its positive cone
representative is \((1_{L_2(M)}\otimes\Ad w)\xi\).
Applying \eqref{eq:AY} to \(v\) and \((b\otimes w)v\) proves
\eqref{eq:rightstandard}.
\end{proof}

\section{A quantitative form of the coercivity argument}\label{sec:quantitative-coercivity}

%We prove Theorem~\ref{thm:gap} using a coercivity estimate on a Walsh
%space. Matrix \(L_2\)-spaces carry normalized traces; the
%Hilbert--Schmidt norm is unnormalized.

\subsection{From a spectral gap estimate to uniform coercivity}\label{sec:quant-abstract}

%If a nonnegative quadratic form controls one orthogonal component, approximate invariance
%under a unitary exchanging the two components also controls the other.
%The following lemma makes this observation quantitative.
We begin with a simple observation that highlights a key feature of the
Walsh space introduced below.

\begin{lemma}\label{lem:quant-complementary}
Let \(\mathcal H\) be a Hilbert space, \(L\in\B(\mathcal H)^+\),
\(P\in\B(\mathcal H)\) an orthogonal projection, and
\(U\in\mathcal U(\mathcal H)\). If
\[
 P\le\gamma L,\qquad U^*PU=1_{\mathcal H}-P
\]
for some \(\gamma>0\), then
\[
 L+(U-1_{\mathcal H})^*(U-1_{\mathcal H})
 \ge(2\gamma+1)^{-1}1_{\mathcal H}.
\]
\end{lemma}
\begin{proof}
For \(\xi\in\mathcal H\), we have
\(\|\xi\|^2=\|P\xi\|^2+\|PU\xi\|^2\). The triangle inequality
in \(\mathcal H\oplus\mathcal H\) therefore gives
\[
 \|\xi\|\le2^{1/2}\|P\xi\|+\|(U-1_{\mathcal H})\xi\|
 \le(2\gamma)^{1/2}\|L^{1/2}\xi\|+\|(U-1_{\mathcal H})\xi\|.
\]
Apply Cauchy--Schwarz.
\end{proof}

Fix \(n\ge2\) and a Hilbert space \(\mathcal K\). For
\(S\subseteq[n]\), let \(p_S\in\mathbb M_n\) be the orthogonal
projection onto the coordinates indexed by \(S\). We consider the Walsh space 
\[
 \mathcal W_n=\bigoplus_{S\subseteq[n]}(\ell_2^n\otimes\mathcal K),
\]
and define the orthogonal projection \(P_n\in\B(\mathcal W_n)\) and the
normalized Walsh number operator \(N_n\in\B(\mathcal W_n)^+\) by
\[
 P_n=\bigoplus_{S\subseteq[n]}(p_S\otimes1_{\mathcal K}),
 \qquad
 N_n=\bigoplus_{S\subseteq[n]}\frac{|S|}{n}\,1_{\ell_2^n\otimes\mathcal K}.
\]
For \(w=(w_S)_{S\subseteq[n]}\in\mathcal W_n\), write
\(w_S=\sum_{k=1}^n e_k\otimes w_{k,S}\), where \((e_k)_{k=1}^n\)
is the canonical basis of \(\ell_2^n\) and \(w_{k,S}\in\mathcal K\).
In these coordinates, for every \(k\in[n]\) and \(S\subseteq[n]\),
\[
 (P_nw)_{k,S}=
 \begin{cases}
  w_{k,S},&k\in S,\\
  0,&k\notin S,
 \end{cases}
 \qquad
 (N_nw)_{k,S}=\frac{|S|}{n}\,w_{k,S}.
\]

The next proposition gives a uniform coercivity estimate on
\(\mathcal W_n\). It combines the matrix spectral gap estimate, the
potential \(N_n\), and a unitary exchanging
\(\operatorname{ran}P_n\) and \(\ker P_n\).

\begin{proposition}\label{prop:quant-walsh-control}
Let \(R_1,R_2\in\mathcal U(\ell_2^n)\) satisfy
\begin{equation}\label{eq:quant-poincare}
 \|x-\tr_n(x)1_n\|_{\mathrm{HS}}^2
 \le C\sum_{i=1}^2\|R_ixR_i^*-x\|_{\mathrm{HS}}^2
 \qquad\text{for every }x\in\mathbb M_n
\end{equation}
for some \(C>0\), where \(\tr_n=n^{-1}\Tr\).
For \(A_1,A_2\in\mathcal U(\mathcal K)\), put
\[
 W_i=\bigoplus_{S\subseteq[n]}(R_i\otimes A_i)
 \in\mathcal U(\mathcal W_n)\qquad(i=1,2).
\]
If \(W_3\in\mathcal U(\mathcal W_n)\) satisfies
\(W_3^*P_nW_3=1_{\mathcal W_n}-P_n\), then
\begin{equation}\label{eq:quant-coercivity}
 \|N_n^{1/2}w\|^2+\sum_{i=1}^3\|W_iw-w\|^2
 \ge(4C+3)^{-1}\|w\|^2
 \qquad\text{for every }w\in\mathcal W_n.
\end{equation}
\end{proposition}
\begin{proof}
Write \(w=(w_S)_{S\subseteq[n]}\in\mathcal W_n\), where
\(w_S\in\ell_2^n\otimes\mathcal K\). Let \(h_S\in\mathbb M_n^+\)
be determined by
\[
 \Tr(xh_S^2)=\langle(x\otimes1_{\mathcal K})w_S,w_S\rangle
 \qquad\text{for every }x\in\mathbb M_n.
\]
Thus \(\|h_S\|_{\mathrm{HS}}=\|w_S\|\). Lemma~\ref{lem:standard},
applied to these vector functionals, gives by homogeneity
\[
 \sum_S\|R_ih_SR_i^*-h_S\|_{\mathrm{HS}}^2
 \le2\|W_iw-w\|^2\qquad(i=1,2).
\]
Consequently, \eqref{eq:quant-poincare} implies
\[
 \sum_S\|h_S-\tr_n(h_S)1_n\|_{\mathrm{HS}}^2
 \le2C\sum_{i=1}^2\|W_iw-w\|^2.
\]
Since
\[
 \|P_nw\|^2=\sum_S\|h_Sp_S\|_{\mathrm{HS}}^2,
 \qquad n\tr_n(h_S)^2\le\|h_S\|_{\mathrm{HS}}^2=\|w_S\|^2,
\]
the triangle inequality yields
\begin{align*}
 \|P_nw\|
 &\le\Bigl(\sum_S\|h_S-\tr_n(h_S)1_n\|_{\mathrm{HS}}^2\Bigr)^{1/2}
       +\Bigl(\sum_S|S|\tr_n(h_S)^2\Bigr)^{1/2}\\
 &\le\Bigl(2C\sum_{i=1}^2\|W_iw-w\|^2\Bigr)^{1/2}
       +\|N_n^{1/2}w\|.
\end{align*}
Put
\[
 L=N_n+\sum_{i=1}^2(W_i-1_{\mathcal W_n})^*(W_i-1_{\mathcal W_n})
 \in\B(\mathcal W_n)^+.
\]
Cauchy--Schwarz gives \(P_n\le(2C+1)L\).
Lemma~\ref{lem:quant-complementary}, applied to \(P_n,L,W_3\),
proves \eqref{eq:quant-coercivity}.
\end{proof}

\subsection{The Walsh space associated with a representation and the triple \((U_{i,n})_{i=1}^3\)}

Fix \(n\ge2\) and a unitary representation
\(\pi:\F_3\times\F_3\to\mathcal U(K)\) on a Hilbert space \(K\).
Put \(M=\pi(\F_3\times\{e\})''\subset\B(K)\). Recall from
Section~2.3 that
\[
 H_n=\ell_2^n\otimes H_n^0\otimes\ell_2^{m_n}.
\]
Define
\[
 \begin{aligned}
 \mathcal L_n&=L_2(M)\otimes L_2(\B(H_n))\\
 &=L_2(M)\otimes L_2(\mathbb M_n)\otimes L_2(\mathbb M_{2^n})
        \otimes L_2(\mathbb M_{m_n}),\\
 \mathcal W_n&=L_2(M)\otimes L_2(\mathbb M_n)
        \otimes L_2(\A_n)\otimes\C1_{m_n}\subset\mathcal L_n.
 \end{aligned}
\]
Here \(1_{m_n}\in L_2(\mathbb M_{m_n})\) has norm one.
Only the first triple \((U_{i,n})_{i=1}^3\) is used in this subsection.
For \(i=1,2,3\), set
\[
 \widetilde W_i=\lambda_M(\pi(u_i,e))\otimes L_{\overline{U_{i,n}}}
 \in\mathcal U(\mathcal L_n),
 \qquad W_i=\widetilde W_i|_{\mathcal W_n}.
\]
The operators \(\lambda_M\) and \(L_x\) are the left actions defined
in Section~\ref{sec:standard}. Bars denote conjugation in the standard
tensor basis of \(H_n\). Formula~\eqref{eq:encoding} shows that
\(\widetilde W_i\mathcal W_n=\mathcal W_n\), so
\(W_i\in\mathcal U(\mathcal W_n)\).

Let \(\Omega_n=\{-1,1\}^n\), with uniform probability measure
\(\mu_n\). Using the standard isometric identifications
\[
 L_2(\mathbb M_n)\simeq\ell_2^n\otimes\ell_2^n,\qquad
 L_2(\A_n)\simeq L^2(\Omega_n,\mu_n),
\]
where \(\chi_{S,n}\) corresponds to
\(\varepsilon\mapsto\prod_{k\in S}\varepsilon_k\), and expanding in the
Walsh basis, we obtain
\[
 \mathcal W_n\simeq\bigoplus_{S\subseteq[n]}
 \bigl(\ell_2^n\otimes(L_2(M)\otimes\ell_2^n)\bigr).
\]
This is the setting of Section~\ref{sec:quant-abstract} with
\(\mathcal K=L_2(M)\otimes\ell_2^n\); we retain the notation \(P_n,N_n\).
On each summand, \(W_i\), for \(i=1,2\), acts as
\[
 \overline{R_{i,n}}\otimes
 \bigl(\lambda_M(\pi(u_i,e))\otimes1_{\ell_2^n}\bigr).
\]
For \(w\in\mathcal W_n\), use the coordinates
\(w_{k,S}\in L_2(M)\otimes\ell_2^n\) of Section~\ref{sec:quant-abstract}.
By \eqref{eq:toggle}, the third unitary satisfies
\[
 (W_3w)_{k,S}
 =\bigl(\lambda_M(\pi(u_3,e))\otimes1_{\ell_2^n}\bigr)
       w_{k,S\triangle\{k\}}
 \qquad(k\in[n],\ S\subseteq[n]).
\]
Thus \(W_3\) exchanges the components with \(k\in S\) and \(k\notin S\),
so \(W_3^*P_nW_3=1_{\mathcal W_n}-P_n\).
Since the matrix gap \eqref{eq:quantumgap} is unchanged by entrywise
conjugation, Proposition~\ref{prop:quant-walsh-control} yields
\begin{equation}\label{eq:quant-concrete-coercivity}
 \|N_n^{1/2}w\|^2+\sum_{i=1}^3\|W_iw-w\|^2
 \ge(4C+3)^{-1}\|w\|^2
 \qquad\text{for every }w\in\mathcal W_n.
\end{equation}

\subsection{From an almost intertwiner to an almost invariant vector in the Walsh space}

The next proposition turns an almost intertwiner into a unit vector
\(w\in\mathcal W_n\) with small displacements under \(W_1,W_2,W_3\) and
small \(\|N_n^{1/2}w\|\). The coercivity estimate \eqref{eq:quant-concrete-coercivity}
then rules out sufficiently small intertwining defects.

\begin{proposition}\label{prop:quant-reduction}
For every \(\varepsilon>0\), there exists \(\eta>0\) such that the
following holds for every \(n\ge2\) and every unitary representation
\(\pi:\F_3\times\F_3\to\mathcal U(K)\) on a Hilbert space \(K\).
If \(T\in\B(H_n,K)\) satisfies
\begin{equation}\label{eq:quant-almost-intertwiners}
 \|T\|_{\mathrm{HS}}=1,\qquad
 \sum_{i=1}^3\left(
 \|\pi(u_i,e)T-TU_{i,n}\|_{\mathrm{HS}}^2
 +\|\pi(e,u_i)T-TV_{i,n}\|_{\mathrm{HS}}^2
 \right)\le\eta,
\end{equation}
then there is a unit vector \(w\in\mathcal W_n\) such that
\begin{equation}\label{eq:quant-small-norms}
 \|N_n^{1/2}w\|^2+\sum_{i=1}^3\|W_iw-w\|^2<\varepsilon.
\end{equation}
\end{proposition}
\begin{proof}
Fix \(\varepsilon>0\) (we will choose \(\eta>0\) at the end). Let us first pass to the standard form.
Identifying isometrically the Hilbert space \(\operatorname{HS}(H_n,K)\)
of Hilbert--Schmidt operators with \(K\otimes\overline{H_n}\), the
almost intertwiner \(T\) gives a unit vector
\(v\in K\otimes\overline{H_n}\) such that, for every \(i=1,2,3\),
\[
 \begin{aligned}
 \|(\pi(u_i,e)\otimes\overline{U_{i,n}})v-v\|
 &=\|\pi(u_i,e)T-TU_{i,n}\|_{\mathrm{HS}},\\
 \|(\pi(e,u_i)\otimes\overline{V_{i,n}})v-v\|
 &=\|\pi(e,u_i)T-TV_{i,n}\|_{\mathrm{HS}}.
 \end{aligned}
\]
Let
\(\xi\in\mathcal L_n\) be the unit vector in the positive cone
representing the vector state of \(v\) on \(M\bar\otimes\B(H_n)\).
Put
\[
 S_i=1_{L_2(M)}\otimes\Ad\overline{V_{i,n}}
 \in\mathcal U(\mathcal L_n)\qquad(i=1,2,3).
\]
Since \(\pi(e,u_i)\in M'\), Lemma~\ref{lem:standard} gives
\begin{equation}\label{eq:quant-standard-actions}
 \sum_{i=1}^3\|\widetilde W_i\xi-\xi\|^2\le\eta,
 \qquad
 \sum_{i=1}^3\|S_i\xi-\xi\|^2\le2\eta.
\end{equation}

Let us now project onto the Walsh space \(\mathcal W_n\).
Let \(h\in\mathcal L_n\) be the orthogonal projection of \(\xi\) onto
\[
 L_2(M)\otimes L_2(\mathbb M_n)\otimes L_2(\mathbb M_{2^n})
       \otimes\C1_{m_n}.
\]
By \eqref{eq:quantumgap}, applied in dimension \(m_n\),
\begin{equation}\label{eq:quant-scalar-projection}
 \|\xi-h\|^2\le C\sum_{i=1}^2\|S_i\xi-\xi\|^2\le2C\eta.
\end{equation}
Consequently,
\[
 \|S_3h-h\|\le(2\eta)^{1/2}+2\|\xi-h\|.
\]
Let \(E_A,E_B\in\B(\mathcal L_n)\) be the orthogonal projections onto
\(\mathcal W_n\) and
\[
 L_2(M)\otimes L_2(\mathbb M_n)\otimes L_2(\mathcal B_n)
       \otimes\C1_{m_n}.
\]
Since \(\overline{V_{3,n}}=V_{3,n}\),
Proposition~\ref{prop:walsh-averaging}, with \(H=L_2(M)\), gives
\begin{equation}\label{eq:quant-two-projections}
 \|h-E_Ah\|^2+\|h-E_Bh\|^2=\|S_3h-h\|^2.
\end{equation}
Set \(w_0=E_Ah=E_A\xi\in\mathcal W_n\).
To apply Proposition~\ref{prop:walsh-degree-control}, observe that
\begin{equation}\label{eq:quant-distance}
 \begin{aligned}
 \operatorname{dist}(w_0,\operatorname{ran}E_B)^2
 &=\|(1_{\mathcal L_n}-E_B)w_0\|^2\\
 &\le2\bigl(\|h-E_Ah\|^2+\|h-E_Bh\|^2\bigr).
 \end{aligned}
\end{equation}
That proposition, with coefficient space
\(L_2(M)\otimes L_2(\mathbb M_n)\), yields
\[
 \|N_n^{1/2}w_0\|^2
 \le(1-e^{-1})^{-1}\operatorname{dist}(w_0,\operatorname{ran}E_B)^2.
\]
Also,
\(\|\xi-w_0\|^2=\|\xi-h\|^2+\|h-E_Ah\|^2\).
The preceding estimates imply, for a constant \(C_0>0\) depending only
on \(C\),
\begin{equation}\label{eq:quant-norm}
 \|\xi-w_0\|^2+\|N_n^{1/2}w_0\|^2\le C_0\eta.
\end{equation}

To conclude, we only need to normalize \(w_0\).
Since \(\widetilde W_i\mathcal W_n=\mathcal W_n\), the projection
\(E_A\) commutes with \(\widetilde W_i\). Thus
\[
 \sum_{i=1}^3\|W_iw_0-w_0\|^2
 =\sum_{i=1}^3\|E_A(\widetilde W_i\xi-\xi)\|^2\le\eta.
\]
Moreover, \(\|w_0\|^2=1-\|\xi-w_0\|^2\ge1-C_0\eta\).
For \(C_0\eta<1\), the unit vector \(w=w_0/\|w_0\|\) satisfies
\[
 \|N_n^{1/2}w\|^2+\sum_{i=1}^3\|W_iw-w\|^2
 \le(1-C_0\eta)^{-1}(C_0+1)\eta<\varepsilon
\]
provided \(\eta>0\) is sufficiently small. This choice is independent
of \(n\), \(K\), and \(\pi\).
\end{proof}

\section{Proofs of the main results}\label{sec:nonllp}

\subsection{Proof of Theorem~\ref{thm:gap}}

Theorem~\ref{thm:gap} follows from \eqref{eq:smallcomm},
\eqref{eq:quant-concrete-coercivity}, and
Proposition~\ref{prop:quant-reduction}.\qed

\subsection{Proof of Theorem~\ref{thm:main}}

We also denote by \(\rho\) the unital \(*\)-homomorphism extending
the representation defined after Theorem~\ref{thm:gap}:
\[
 \rho:C^*(\F_3\times\F_3)\longrightarrow
 \ell_\infty(\B(H_n))/c_0(\B(H_n)).
\]
Suppose that \(C^*(\F_3\times\F_3)\) has the LLP. Lift \(\rho\) on
the operator system spanned by \((u_i,e)\) and \((e,u_i)\), \(i=1,2,3\),
and extend the coordinate maps by Arveson's theorem. This gives u.c.p.\
maps \(\beta_n:C^*(\F_3\times\F_3)\to\B(H_n)\) satisfying
\[
 \|\beta_n(u_i,e)-U_{i,n}\|\longrightarrow0,\qquad
 \|\beta_n(e,u_i)-V_{i,n}\|\longrightarrow0
 \qquad(i=1,2,3).
\]
Choose Stinespring dilations \(\beta_n(x)=J_n^*\pi_n(x)J_n\), where
\(\pi_n:C^*(\F_3\times\F_3)\to\B(K_n)\) is a unital
\(*\)-representation and \(J_n\in\B(H_n,K_n)\) is an isometry.
For \(g\in\F_3\times\F_3\) and \(X\in\mathcal U(H_n)\),
\begin{equation}\label{eq:stinespring-estimate}
 (\dim H_n)^{-1}\|\pi_n(g)J_n-J_nX\|_{\mathrm{HS}}^2
 =2\Ree\tr_{\dim H_n}(1-X^*\beta_n(g))
 \le2\|\beta_n(g)-X\|_2.
\end{equation}
Here \(\|\cdot\|_2\) is the normalized Hilbert--Schmidt norm.
Apply \eqref{eq:gap} with \(T=J_n\). Since
\(\|J_n\|_{\mathrm{HS}}^2=\dim H_n\) and \(\|\cdot\|_2\le\|\cdot\|\),
\eqref{eq:stinespring-estimate} gives
\[
 c\le2\sum_{i=1}^3\bigl(
 \|\beta_n(u_i,e)-U_{i,n}\|+\|\beta_n(e,u_i)-V_{i,n}\|\bigr)
 \longrightarrow0,
\]
a contradiction.

The inclusions
\(\F_3\times\F_3\subset\F_2\times\F_2\subset\F_k\times\F_m\),
for \(k,m\ge2\), follow from the index-two embedding
\(\F_3\hookrightarrow\F_2\) and the first two generators of each
free group. Since the LLP passes to full group \(C^*\)-algebras of
subgroups \cite[Corollary~3.17]{FFW}, the theorem follows.\qed

\subsection{Proof of Theorem~\ref{thm:local-embedding}}
\label{sec:local-embedding}

Recall that a tracial state \(\tr_A\) on a unital \(C^*\)-algebra \(A\)
is \emph{amenable} if its GNS representation
\(\pi_{\tr_A}:A\to\pi_{\tr_A}(A)''\) extends to a u.c.p.\ map on
\(\B(K)\), for a faithful unital representation \(A\subset\B(K)\)
\cite[Theorem~3.1.6(5)]{Brown}. We use the following consequence of
local embeddability from \cite{Roydor}.

\begin{proposition}\label{prop:tracial-local-lifting}
Let \(A\) be a unital \(C^*\)-algebra that locally embeds into
\(C^*(\F_\infty)\). Let \(\alpha:A\to M\) be a unital
\(*\)-homomorphism, where \(M\) is a finite von Neumann algebra with
a faithful normal tracial state \(\tr_M\), and suppose that
\(\tr_M\alpha\) is amenable. For every unital \(C^*\)-algebra \(D\),
unital surjective \(*\)-homomorphism \(q:D\to M\), and
finite-dimensional operator system \(E\subset A\), there is a u.c.p.\
map \(\beta:E\to D\) with \(q\beta=\alpha|_E\).
\end{proposition}
\begin{proof}
Represent \(A\) faithfully in \(\B(K)\), and put
\(\tr_A=\tr_M\alpha\). The trace-preserving identification
\(\pi_{\tr_A}(A)''\cong\alpha(A)''\subset M\) and amenability give a
u.c.p.\ extension \(\widetilde\alpha:\B(K)\to M\) of \(\alpha\).
The local lifting characterization of local embeddability
\cite[Theorem~20.5(v)]{Pisier}, applied to
\(A\subset\B(K)\xrightarrow{\widetilde\alpha}M\cong D/\ker q\), gives local
completely contractive liftings. Since \(\alpha\) is u.c.p.,
\cite[Corollary~9.47]{Pisier} gives a u.c.p.\ lifting on \(E\).
\end{proof}

\begin{proof}[Proof of Theorem~\ref{thm:local-embedding}]
Residual finiteness provides finite quotients \(\Lambda_n\) of
\(\F_3\times\F_3\) and associated regular representations
\(r_n:\F_3\times\F_3\to\mathcal U(\ell_2(\Lambda_n))\) such that
\(r_n(g)\) has trace zero eventually for each \(g\ne(e,e)\).
On \(\widehat H_n=H_n\otimes\ell_2(\Lambda_n)\), put
\[
 \widehat U_{i,n}=U_{i,n}\otimes r_n(u_i,e),\qquad
 \widehat V_{i,n}=V_{i,n}\otimes r_n(e,u_i)
 \qquad(i=1,2,3).
\]
These unitaries still have commutators tending to zero in operator norm
and satisfy \eqref{eq:gap} with the same constant \(c\).
To verify the latter, fix a representation
\(\pi:\F_3\times\F_3\to\mathcal U(K)\) and
\(T\in\B(\widehat H_n,K)\). Define
\[
 \Theta_n(T)\xi=\sum_{h\in \Lambda_n}T(\xi\otimes\delta_h)
       \otimes\overline{\delta_h}\qquad(\xi\in H_n),
\]
where \((\delta_h)_{h\in \Lambda_n}\) is the canonical basis of \(\ell_2(\Lambda_n)\).
Then \(\Theta_n(T)\in\B(H_n,K\otimes\overline{\ell_2(\Lambda_n)})\),
\(\|\Theta_n(T)\|_{\mathrm{HS}}=\|T\|_{\mathrm{HS}}\), and, for every
\(g\in\F_3\times\F_3\) and \(X\in\B(H_n)\),
\[
 \begin{aligned}
 &\|(\pi(g)\otimes\overline{r_n(g)})\Theta_n(T)-\Theta_n(T)X\|_{\mathrm{HS}}\\
 &\hspace{20mm}=\|\pi(g)T-T(X\otimes r_n(g))\|_{\mathrm{HS}}.
 \end{aligned}
\]
Apply \eqref{eq:gap} to \(\pi\otimes\overline{r_n}\) and
\(\Theta_n(T)\).

Fix a free ultrafilter \(\omega\) on \(\mathbb N\), and let
\(\|\cdot\|_{2,n}\) denote the normalized Hilbert--Schmidt norm on
\(\B(\widehat H_n)\). Define the tracial ultraproduct
\[
 \begin{aligned}
 M_\omega&=\ell_\infty(\B(\widehat H_n))/I_\omega,\\
 I_\omega&=\{(X_n)\in\ell_\infty(\B(\widehat H_n)):
                  \lim_{n\to\omega}\|X_n\|_{2,n}=0\},
 \end{aligned}
\]
with quotient map \(q_\omega\) and faithful normal trace
\[
 \tr_\omega(q_\omega((X_n)))
 =\lim_{n\to\omega}(\dim\widehat H_n)^{-1}\Tr(X_n).
\]
The amplified matrices define a unital \(*\)-homomorphism
\(\alpha:C^*(\F_3\times\F_3)\to M_\omega\) by
\[
 \alpha(u_i,e)=q_\omega((\widehat U_{i,n})_n),\qquad
 \alpha(e,u_i)=q_\omega((\widehat V_{i,n})_n)
 \quad(i=1,2,3).
\]
Evaluating a word in these matrices shows that
\(\tr_\omega\alpha=\tr_{\mathrm{can}}\), where \(\tr_{\mathrm{can}}\) is the canonical trace:
\(\tr_{\mathrm{can}}(e,e)=1\) and \(\tr_{\mathrm{can}}(g)=0\) for \(g\ne(e,e)\).
The normalized traces of the finite-dimensional representations \(r_n\)
converge pointwise to \(\tr_{\mathrm{can}}\), so \(\tr_{\mathrm{can}}\) is amenable
\cite[Theorem~3.1.6(2)]{Brown}.

Suppose that \(C^*(\F_3\times\F_3)\) locally embeds into
\(C^*(\F_\infty)\). Proposition~\ref{prop:tracial-local-lifting}
lifts \(\alpha\) through \(q_\omega\) on the operator system spanned
by the six generators. Arveson's theorem extends the coordinate maps to
u.c.p.\ maps \(\beta_n:C^*(\F_3\times\F_3)\to\B(\widehat H_n)\), with
\[
 \lim_{n\to\omega}\|\beta_n(u_i,e)-\widehat U_{i,n}\|_{2,n}=0,
 \qquad
 \lim_{n\to\omega}\|\beta_n(e,u_i)-\widehat V_{i,n}\|_{2,n}=0
 \quad(i=1,2,3).
\]
Apply the amplified estimate \eqref{eq:gap} to Stinespring dilations of
\(\beta_n\). The calculation \eqref{eq:stinespring-estimate} gives
\[
 c\le2\sum_{i=1}^3\bigl(
 \|\beta_n(u_i,e)-\widehat U_{i,n}\|_{2,n}
 +\|\beta_n(e,u_i)-\widehat V_{i,n}\|_{2,n}\bigr)
 \xrightarrow[n\to\omega]{}0,
\]
a contradiction. Finally, local embeddability passes to operator
subspaces, and the inclusions
\(\F_3\times\F_3\subset\F_2\times\F_2\subset\F_k\times\F_m\)
induce isometric \(*\)-embeddings of full group \(C^*\)-algebras.
\end{proof}

\begin{remark}[Further groups]
It is not very interesting to make an exhaustive catalogue of groups to
which our results apply, but we record a few examples.
If \(G\) contains \(\F_2\times\F_2\), then \(C^*(G)\) fails the LLP
and does not locally embed into \(C^*(\F_\infty)\): the LLP passes to
full group \(C^*\)-algebras of subgroups \cite[Corollary~3.17]{FFW},
and local embeddability passes to operator subspaces.
Examples include products \(G_1\times G_2\) with \(\F_2\subset G_i\);
\(\SL_n(\Z)\) for \(n\ge4\) and \(\operatorname{Sp}_{2n}(\Z)\)
for \(n\ge2\); the braid and pure braid groups \(B_n,PB_n\) for
\(n\ge5\) \cite[Lemma~12]{KP}; the mapping class groups of closed
oriented surfaces of genus at least two; \(\operatorname{Aut}(\F_n)\)
and \(\operatorname{Out}(\F_n)\) for \(n\ge3\)
\cite[Section~1]{BW}; and right-angled Artin groups whose defining
graphs contain an induced cycle of length four. The same conclusions
hold for their finite-index subgroups.
For \(\SL_3(\Z)\), they follow instead from
\cite[Corollary~B]{ISW} and \cite{Roydor}, respectively.
\end{remark}

\medskip

\begin{samepage}
\noindent \textbf{Acknowledgments.} We thank Gilles Pisier for many stimulating discussions on questions concerning the WEP and the LLP.\\
\textbf{AI Assistance.} ChatGPT was used to develop certain proofs through an
interactive dialogue. We did not use Lean to check the proofs. A year ago,
we tried to adapt Ioana's work on Hilbert--Schmidt stability to this LLP
question. We thought that we should retain the spectral gap estimates
while changing the deformation, but we could not make this idea work.
A few months ago, we returned to the question with \cite{Roydor} in mind
and spent a few weeks trying to prove that \(C^*(\F_2\times\F_2)\)
locally embeds into \(C^*(\F_\infty)\), that is,
\(d_f(C^*(\F_2\times\F_2))=1\) in Pisier's notation
\cite[Section~20]{Pisier}. This is equivalent to \(d_f(S)=1\), where
\(S\) is the operator system spanned by the four canonical unitary
generators. One can prove that \(d_f(S)<1.072\), which may help explain
the difficulty of the LLP question. We then returned to our original
idea and suggested that ChatGPT explore other deformations or semigroups.
We eventually suggested using the Walsh semigroup (see, for example,
Lust-Piquard's work) to address the question. The idea of proving
Theorem~\ref{thm:local-embedding} is taken from \cite{Roydor}.\par
\end{samepage}

\end{document}